\documentclass[12pt]{article}
\usepackage{pudlatex}
\usepackage{hyperref} 

\newcommand\tw{\mathop{\rm tw}\nolimits}

\title{Colorings of $k$-sets with low discrepancy on small sets}

\author{Pavel Pudl\'ak
\thanks{partially supported by grants 19-27871X and 23-048255 of the Czech Grant Agency and the institute grant RVO: 67985840.} \\
Institute of Mathematics, CAS, Prague, Czech Republic
\and 
Vojt\v{e}ch R\"odl
\thanks{partially supported by NSF grants DMS 1764385 and DMS 2300347.}\\
Emory University, Atlanta, USA
}

\begin{document}

\maketitle

\begin{abstract}
For $0<\delta\leq 1$, let $R_k(m;\delta)$ denote the smallest $N$ such that every coloring of $k$-element subsets by two colors yields an $m$-element set $M$ with relative discrepancy $\delta$, which means that one color class has at least $(\frac{1+\delta}2){m\choose k}$ elements. The number $R_k(m;\delta)$ may be viewed as an extension of the usual $k$-hypergraph Ramsey number because  $R_k(m)=R_k(m,1)$. Our main result is the following theorem.
{\it For some constants $c,k_0$, and $\eps>0$, and for all $k\geq k_0$, $c\log k\leq n\leq k/11$,}
\[
R_k(k+n;2^{-\eps n})\geq \tw_{\lfloor k/n\rfloor}(2).
\]
In particular, for $n=\lceil c\log k\rceil$, we get a tower of height
$\delta k/\log k$ and relative discrepancy polynomial in~$k$.

\end{abstract}

\section{Introduction}

A well-known theorem of F. P. Ramsey~\cite{ramsey} ensures for every positive integers $k<m$ the existence of an integer $N$ with the property that for any coloring of all $k$-element subsets of a set $[N]=\{1,2\dts N\}$ by two colors, there exists an $m$-element set $M\sub[N]$ with all ${m\choose k}$ $k$-element subsets colored by the same color. The Ramsey number $R_k(m)$ is then defined as the least integer $N$ with the above property. The numbers $R_k(m)$ have been studied for almost ninety years. While for $k=2$, it is known that $R_2(n)$ grows exponentially, there is still an exponential gap between the lower and upper bounds for $k\geq 3$, and the best bounds are
\bel{e-*}
tw_{k-1}(c_1m^2)\leq R_k(m)\leq tw_k(c_2m),
\ee
where $c_1,c_2$ are absolute constants and $tw_k(x)$ is defined by $tw_1(x)=x$ and $tw_{j+1}=2^{tw_j(x)}$, see~\cite{erdos-hajnal,erdos-rado,erdos-szekeres,graham-rothchild}.

Given $k<m$ and $t$ with $\frac 12{m\choose k}<t\leq{m\choose k}$, a natural extension of Ramsey numbers $R_k(m)$ is, for every 2-coloring of ${[N]\choose k}$, to require the existence of a set $M$, $|M|=m$, inducing at least $t$ monochromatic $k$-tuples. To state the results about this generalization, we will use the following notation.
For $0<\delta\leq 1$, $R_k(m;\delta)$ denotes the smallest $N$ such that every coloring of $k$-element subsets by two colors yields an $m$-element set $M$ in which one color class has at least $(\frac{1+\delta}2){m\choose k}$ elements. In other words, $\delta$ is a bound on the relative discrepancy.%
\footnote{In~\cite{erdos} Erd\H os introduced a function $F_2^{(k)}(N,\alpha)$, for $0\leq\alpha\leq 1/2$, which is in a sense, inverse to our $R_k(m;\delta)$. This function is defined as the largest $m$ such that every 2-coloring of $[N]$ has a set whose relative discrepancy is $1-2\alpha$.}
Clearly,  $R_k(m)=R_k(m,1)$. Our main result is the following theorem.

\bt\label{t-1}
There exist constants $c,k_0$, and $\eps$ such that for all $k\geq k_0$, $c\log k\leq n\leq k/4$,
\bel{e-main}
R_k(k+n;2^{-\eps n})\geq \tw_{\lfloor k/n\rfloor}(2).
\ee
\et
In particular, for $n=\lceil c\log k\rceil$, we get a tower of height $\delta k/\log k$, $\delta>0$, and polynomial discrepancy.

Questions concerning discrepancy in Ramsey theory have been studied before. Let us mention a well-known result of P.~Erd\H os and J.~Spencer \cite{erdos}, which states that for a given $k$ and $0<\delta<\delta(k)$ sufficiently small,
\bel{E-S}
2^{c_1m^{k-1}}\leq R_k(m;\delta)\leq 2^{c_2m^{k-1}},
\ee
where $c_1,c_2$ depend on $\delta$ and $k$. Since $R_k(m,1)$ is the usual Ramsey number for which the bounds are of tower type, Erd\H os asked whether the transition from $\delta$ close to $0$ to $\delta=1$ is continuous, or there are jumps (see~\cite{chung-graham}). This problem is still open in full generality, though some progress has been made. In particular, D.~Conlon, J.~Fox, and B. Sudakov proved that for every $\delta<1$, there exists $c$ such that $R_3(m;\delta)\leq 2^{cm^{2}}$~\cite{conlon-fox-sudakov-Israel}. 
Our result only shows that  the constants $c_1,c_2$ in (\ref{E-S}) must grow at least as a tower of twos of height ${\eps' k/\log k}$. Since our lower bound is a function of $k$ and not of $m$, it says nothing about how $R_k(m;\delta)$ depends on~$m$.

To prove a lower bound on Ramsey numbers, one has to exhibit a construction of a coloring. The coloring that we construct in the proof of our main theorem is not quite explicit because it is composed of two mappings, one of which is a random function. Therefore, we also present an explicit construction, which, however, yields a weaker lower bound and only works for an odd number of colors.

The main tools for proving lower bounds on $k$-uniform hypergraph
Ramsey numbers have been either the probabilistic method or the
stepping-up lemma of Erd\H os and Hajnal~\cite{erdos-hajnal}.
In the regime when $m-k$ is very small, a method based on shift-graphs~\cite{duffus-lefmann-rodl}, or the approach taken by G.~Moshkovitz and A.~Shapira~\cite{moshkovitz} seems to yield better bounds. 
In this paper, we use colorings of shift-graphs; we do not know how to prove such a result using the stepping-up lemma.\footnote{In \cite{pudlak-rodl} we used the stepping-up lemma to get triple exponential from a double exponential for a related problem, but we could not iterate this construction to get more.}

This paper builds on ideas from our previous paper~\cite{pudlak-rodl}, where we constructed extractors of a particular type. Studying discrepancy from the point of view of the theory of extractors turned out to be very useful.

\paragraph{Acknowledgment.} We are grateful to the referees for their suggestions on improving the presentation of the results.


\section{Preliminaries}

For a positive integer $n$, we denote by $[n]:=\{1,2\dts n\}$. For a positive integer $k$ and a set $S$ of size at least $k$, we denote by ${S\choose k}$ the set of all $k$-element subsets of $S$.
We will abbreviate \emph{``$k$-element sets''} with \emph{``$k$-sets''}.
We will use the convention that when denoting a set of numbers by $\{x_1,x_2\dts x_t\}$, we will always tacitly assume that $x_1<x_2<\dots< x_t$.

\bdf\label{d-1}
Let $f:A\to C$ and $|C|>1$. We will call $f$ a $C$-coloring of $A$ and $C$ the set of colors. The \emph{discrepancy of $f$} is defined by
\[
\max_{i\in C}|\ |f^{-1}(i)|-|A|/|C|\ |.
\]
Note that if $C=\{-1,1\}$, then the discrepancy is $\frac 12|\sum_{x\in A} f(x)|$.

The \emph{relative discrepancy of $f$} is the discrepancy normalized by $|A|$. In terms of probabilities, the relative discrepancy of $f$ is
\[
disc(f,A):=\max_{i\in[c]}|\prob[f(x)=i]-\mbox{$\frac 1{|C|}$}|,
\]
where the probability is w.r.t. the uniform distribution of $x$ on $A$.
\edf

We will use the following simple facts about relative discrepancy. 
\bfa\label{fact}
\ben
\item If $f$ has relative discrepancy $\leq \eps$ on every $A\sub U$, $|A|=m$, then $f$ has relative discrepancy $\leq \eps$ on every $B\sub U$, $|B|\geq m$.
\item Let $A$ be a disjoint union of sets $A_i$, and suppose that the relative discrepancy of $f$ on each $A_i$ is $\leq\eps$. Then the relative discrepancy of $f$ on the entire $A$ is also $\leq\eps$.
\item\label{it-3} Let $B\sub A$ and the relative discrepancy of $f$ on $B$ is $\leq\eps$. Then the relative discrepancy of $f$ on $A$ is $\leq \eps+\frac{|A\setminus B|}{|A|}$. 
\een
\efa

\bdf
Let $l<N$. The \emph{shift graph} $Sh(N,l)$ is the graph on ${[N]\choose l}$ whose edges are pairs of $l$-sets in a \emph{shift position}, i.e., pairs of the form
\[
(\{x_1,x_2\dts x_l\},\{x_2\dts x_l,x_{l+1}\})
\]
where $x_1<x_2<\dots<x_l<x_{l+1}$.
\edf

\bll{l-shift}
Let $l\leq N\leq 4\tw_{l-4}(2)$. Then the chromatic number of $Sh(N,l)$ is~$\leq 3$ and, for $N\geq 2l+1$, equality holds. Moreover, such a $3$-coloring $\chi$ can be explicitly defined.
\el
We postpone the proof to Section~\ref{s-shift}.

We will use some well-known concentration bounds. The following are two inequalities from Theorem~2.3 of~\cite{mcdiarmind} slightly restated.

\bt[Chernoff's bounds]\label{t-chernoff}
Let $X$ be the binomial distribution $Bi(n,p)$. Then
\ben
\item[(a)] $\prob[X\leq \mathbb{E}[X]-\lambda]\leq\exp(-\frac{\lambda^2}{2\mathbb{E}[X]})$;
\item[(b)] if $X=Bi(n,1/2)$, then
$\prob|[X-\mathbb{E}[X]|\geq \lambda]\leq 2\exp(-\frac{2\lambda^2}{n})$.
\een
\et

Furthermore, we will need Hoeffding's bound on the tail of hypergeometric distribution~\cite{hoeffding,serfling}. Since we will  apply it for $k>N/2$, we will state it in the form that is better suited for such a case.

\bt\label{t-hyper}
Let $X$ be the hypergeometric distribution with population size $L$, the number of success states $K$, and the number of draws~$k$.\footnote{Thus $X$ is the random variable counting the number of success states when $k$ elements are drawn.} Let $p:=K/L$ and let $0<\tau<p$. Then
\bel{e-hyper}
\prob[X\leq (p-\tau)k]\leq \exp(-2\tau^2\frac{k^2}{L-k}).
\ee
\et

Next theorem is an application of Azuma's inequality to a Doob martingale. For the Doob's martingale process and Azuma's inequality, see~\cite{mcdiarmind}, page 221.
\bt[Azuma's inequality]\label{l-azuma}
Let $f(x_1\dts x_n)$ be an arbitrary real function. Let $X_1\dts X_n$ be arbitrary random variables. Let $K$ be a constant such that for all $i=1\dts n$ and all $x_1\dts x_n$
\bel{e-azuma}
| \mathbb{E}[f(x_1\dts x_{i-1},x_i,X_{i+1}\dts X_n)]-  \mathbb{E}[f(x_1\dts x_{i-1},X_i,X_{i+1}\dts X_n)]|\leq K,
\ee
where the first expectation is over $X_{i+1}\dts X_n$ and the second over $X_{i}\dts X_n$. Then
\bel{e-azuma2}
\prob[\ |f(X_1\dts X_n)-\mathbb{E}[f(X_1\dts X_n)]|>\lambda\ ]<
2\exp(-\lambda^2/2K^2n),
\ee
where the expectation is over all $X_i$s.
\et

\noindent The bound (\ref{e-azuma}), clearly, follows from 
\bel{e-m1}
| \mathbb{E}[f(x_1\dts x_{i-1},x_i,X_{i+1}\dts X_n)]-  \mathbb{E}[f(x_1\dts x_{i-1},x'_i,X_{i+1}\dts X_n)]|\leq K,
\ee
and it is easier to use the latter one whenever possible.


\section{Proof of Theorem \ref{t-1}}


\subsection{A high-level idea of the proof}

First we prove an interim statement under the additional assumption of
$k$ being an odd multiple of~$n$. 
Proving a lower bound means to construct a coloring with low discrepancy on all sets of size~$k+n$. Thus, for the sake of clarity, we state the lower bound as a statement about coloring $k$-sets.

\bt\label{t-x1}
For every $k$ sufficiently large and for $l$ such that $5\leq l< k/(165\log 2k)$, if $k$ is divisible by $2l-1$, then the following holds true. 
If
\[
 N\leq4\tw_{l-4}(2),
\]
then there exists a coloring $\gamma:{[N]\choose k}\to\{-1,1\}$ such that for every subset $S\sub [N]$ of size $|S|\geq k(1+\frac 1{2l-1})$,
\[
disc(\gamma,\mbox{${S\choose k})$}\leq 2^{-k/213l}.
\]
\et
In terms of the function $R_k(m,\epsilon)$ this theorem says that for  $k$ sufficiently large, $l$ such that $5\leq l< k/(165\log 2k)$, and $k$  divisible by $2l-1$,
\[
R_k(k(1+\frac 1{2l-1}),2^{-k/213l})>4\tw_{l-4}(2).
\]
(The constants are surely not the best that one can get using our proof.)

Let $S\sub [N]$,  $|S|= k(1+\frac 1{2l-1})$.
In the proof of Theorem~\ref{t-x1} we will show that ${S\choose k}$ can be decomposed into a disjoint union of certain sets, which we call \emph{cubes}. Then we show that every element, except of an exponentially small fraction of elements of ${S\choose k}$, is in a sufficiently large cube. Consequently, it suffices to construct a coloring that has low discrepancy on every sufficiently large cube.

Given $N$ and $k$, we construct our coloring in two steps. First we take $l\ll k$ and a coloring $\kappa$ of the shift graph $Sh(N,l)$. If $N$ is not too large, then we can find such a coloring with only three colors (see Lemma~\ref{l-shift}). For a given element $\{x_1\dts x_k\}\in{[N]\choose k}$, where $x_i$ are in the increasing order, we take the string
\[
(\kappa(x_1\dts x_l),\kappa(x_2\dts x_{l+1})\dts \kappa(x_{k-l+1}\dts x_k)).
\]
This defines a mapping $\phi:{[N]\choose k}\to [3]^{k-l+1}$. The cubes will be  defined in such a way that for every two different elements in a cube
\[
\{x_1,x_2\dts x_k\},\{y_1,y_2\dts y_k\},
\]
there exists $i$, $1\leq i\leq k-l+1$, such that the blocks of length $l$
\[
\{x_i,x_{i+1}\dts x_{i+l-1}\},\{y_i,y_{i+1}\dts y_{i+l-1}\}
\]
are in a shift position and hence they have different colors. This implies that $\phi$ is one-to-one on every cube.

Since $\phi$ is one-to-one, the image $\phi(C)$ of a large cube is also large, hence for a random $\psi:[3]^{k-l+1}\to[2]$, $\psi(\phi(C))$ has small discrepancy (i.e., the preimages of 0 and 1 have''almost the same'' size) with high probability. We would like to use the union bound to show that random $\psi$ has low discrepancy on all images of large cubes, but the mere fact that the image $\phi(C)$ is large does not suffice. 
For that we bound the number of sets that could be $\phi$-images of cubes. We are able to get a good bound, because cubes have a very special structure that is partially inherited when we map them into~$[3]^{k-l+1}$.

In Section~\ref{s-explicit} we present a slightly different construction of the coloring that is more explicit. Instead of the random $\psi$, we just add colors modulo an odd number. Thus we only get colorings with odd numbers of colors and furthermore we get a worse bound on how large $N$ can be with respect to $k$. In that construction we do not use all $k-l+1$ blocks of consecutive elements of a $k$-set $X$, instead we divide $X$ into $k/l$ disjoint blocks of size $l$. We also use a different kind of cubes.

We will now proceed to prove Theorem \ref{t-x1}.



\subsection{Covering ${S\choose k}$ by cubes}

\paragraph{Set-up of parameters.} Let $k,l\in\N$ be given, $k$ divisible by $2l-1$. We set
\begin{align*}
n:=& k/(2l-1)\\
m:=& 2ln=(1+\frac 1{2l-1})k.
\end{align*}
For a set $S\sub [N]$, $|S|=m$, we partition $S$ into $n$ consecutive intervals of length $2l$,
\[
S=S_1\cup S_2\cup\dots\cup S_n,
\]
so $\max S_i<\min S_{i+1}$, for $i=1\dts n$.%
\footnote{We call $S_i$ \emph{intervals} to distinguish them from \emph{blocks} in $k$-sets,
although $S_i$s are not intervals in $\N$, but rather intervals in the ordered set $S$.}

\bdf
For a subset $X\sub S$, we will say that \emph{$X$ properly hits $S_i$}, if $|X\cap S_i|=2l-1$. 
\edf
The idea behind this definition is as follows. If $X$ properly hits $S_i$, then there exists exactly one block $X'$ of consecutive elements of $X$ of length~$l$ contained in $S_i$ such that one end of $X'$ is marked with the missing element, which determines $X'$.

We will show that a positive fraction of intervals $S_i$ is properly hit by random $k$-sets with high probability.

\bll{l-number-of-hit}
Let $X$ be a random $k$-subset of $S$. 
Let $Z(X)$ be the number of intervals properly hit by $X$.
\bel{e-2}
\prob[Z(X)<n/2\e ] <\exp(-n/30). 
\ee
\el
\bprf
As it is easier to work with a binomial distribution, we will first derive a bound for sets $X\sub S$ in the binomial distribution.
Our binomial distribution is defined by taking elements randomly independently with probability
\[
p:=\prob[x\in Y]=k/m=1-1/2l.
\]
The probability of properly hitting a single interval for a set $Y$ sampled from this distribution is
\[
\prob[Y\mbox{ properly hits }S_i]=
2l\cdot\frac 1{2l}\cdot\left(1-\frac 1{2l}\right)^{2l-1}\geq 1/\e.
\]
Thus the expected value of $Z(Y)$ is~
$n/\e$. 
Using the Chernoff bound (Theorem~\ref{t-chernoff})
we get
\begin{align*}
  \prob[\ Z(Y) <\mbox{$\frac 12$}\mathbb{E}[Z]\ ]
  &\leq \exp(\mbox{$-\frac 12\cdot\frac 14$}\mathbb{E}[Z])\\
  &\leq\exp(-n/8\e).
\end{align*}
We have chosen the probability distribution so that  $|Y|=k$ is the most likely event among all $|Y|=i$. Indeed, it is well-known that the mode (the quantity we are interested in) of the binomial distribution is either $\lfloor(m+1)p\rfloor$ or $\lceil(m+1)p\rceil-1$. Since $pm=k$ is an integer, the mode is~$k$. Hence $\prob[|Y|=k]\geq\frac 1{m+1}$.%
\footnote{While this bound is not optimal, improving it would have negligible influence on our estimates.}
Thus for the distribution with all sets $Y$ having size~$k$, we get
\begin{align*}
\prob[Z(Y)<n/2\e\ |\ |Y|=k]&=\prob[Z(Y)<n/2\e\ \wedge\ |Y|=k]\cdot\prob[|Y|=k]^{-1}\\
& < \exp(-\frac n{8\e})\cdot(m+1)\\
& = \exp\left(-\frac n{8\e}+\ln(m+1)\right).\qquad (*)
\end{align*}
The expression $(*)$ can be further bounded from above by
\[
 \exp(-\frac n{22}+\ln 2k)\leq \exp(-\frac n{30}),
\]
since $m+1\leq 2k$,
$8\e<22$, and
\[
-\frac n{22}+\frac n{30}+\ln 2k=
-\frac 2{165}\frac k{2l-1}+\ln 2k<
-\frac k{165l}+\ln 2k<0,
\]
where the last inequality is an assumption of Theorem~\ref{t-x1}.
\eprf

\bdf\label{d-cube}
Let $J\sub[n]$ be nonempty and $R\sub S$ such that
\bi
\item for every $j\in J$, $|R\cap S_j|=2l-2$ and the gap between the two missing elements is $l-1$,
\item $|R|+|J|=k$.
\ei
Then \emph{the cube $C(J,R)$} is the set of $k$-element subsets $X\sub S$ such that
\bi
\item for every $j\in J$, $|X\cap S_j|=2l-1$ (i.e., $X$ properly hits $S_j$), and the missing element in $S_j$ is one of the missing elements of $R$ in $S_j$;  
\item for all $j\in[n]\setminus J$, $X\cap S_j=R\cap S_j$.
\ei
The \emph{dimension of the cube} is $|J|$ and the number of elements is, clearly, $2^{|J|}$.
\edf

In other words, for given $R\sub S$ satisfying conditions of Definition~\ref{d-cube}, the elements of $C(J,R)$ are obtained by adding one element to $R$ in every $S_j$ for every $j\in J$.



\bll{l-in-cube}
Every $k$-set $X$ that properly hits at least one $S_i$ is contained in a cube. The maximal cube in which $X$ is contained is unique and has the dimension equal to the number of intervals $S_i$ properly hit by $X$. Hence maximal cubes are disjoint.
\el
\bprf
For a $k$-set $X$ that properly hits at least one $S_i$, let $J$ be the set of the indices of intervals $S_i$ properly hit. We define $R$ as follows.
\bi
\item For $i\not\in J$, put $R\cap S_i:=X\cap S_i$.
\item Now suppose $i\in J$. Let $S_i=\{a_1\dts a_{2l}\}$ and let $a_j$ be the unique element in $S_i\setminus X$. We define
\ben
\item $R\cap S_i:=S_i\setminus\{a_j,a_{j+l}\}$, if $j\leq l$, and
\item $R\cap S_i:=S_i\setminus\{a_{j-l},a_j\}$, if $j>l$.
\een
\ei
$R$ satisfies the condition of the definition and $X$ is in $C(J,R)$. Moreover, $C(J,R)$ is maximal. $C(J,R)$ is a unique maximal cube containing~$X$, because $J$ cannot be larger and in every cube that contains $X$, $R\cap S_i$ has to be as defined above.
\eprf

\begin{corollary}\label{c3.1}
For $k$ sufficiently large,
all sets in ${S\choose k}$, except at most a fraction
$\exp(-\frac n{30})$
of them, can be covered by disjoint cubes of dimensions at least $n/2\e$. \end{corollary}



\subsection{Definition of the coloring}

We will construct our coloring as a composition of two colorings
\[
\phi:{[N]\choose l}\to [3]^{p}\quad\mbox{ and }\quad \psi:[3]^p\to[2],
\]
where $p=k-l+1$.


We will now suppose that $N$ and $l$ satisfy the condition of Lemma~\ref{l-shift} and denote by $\kappa$ a $3$-coloring of $Sh(N,l)$. 
For $X\in{[N]\choose k}$, $X=\{x_1,x_2\dts x_k\}$, and $j\leq k-l+1$, we define \emph{the $j$-th projection} by
\[
\Pi_j(X):=\{x_j,x_{j+1}\dts x_{j+l-1}\},
\]
and we will abbreviate $\Pi_j(X)$ with $X_j$ when there is no danger of confusion.  Our first coloring $\phi$
is defined by
\[
\phi(X):=(\kappa(X_1),\kappa(X_2)\dts\kappa(X_p)).
\]
Note that $X_1\dts X_p$ are all segments of $l$ consecutive elements of $X$.

Our next goal will be to construct a mapping
\[
\psi:[3]^p\to[2]
\]
so that the coloring obtained by composition of the two
\[
\gamma:=\psi\circ\phi:{[N]\choose k}\to [2]
\]
has the properties promised in Theorem~\ref{t-x1}. Towards this end, we need to show that
\ben
\item $\phi$ restricted to each cube $C(J,R)$ is a bijection, and  
\item  due to the special structure of cubes $C(J,R)\sub {S\choose k}$,  $S\in{[N]\choose m}$, the number of $\phi$-images of cubes in $[3]^p$ is bounded by an exponential function in~$n$.
\een

\subsection{Images of cubes $C(J,R)$.}

For $S\in{[N]\choose m}$, $R\sub S$, and $\emptyset\neq J\sub [n]$, consider the cube $C(J,R)$. Let $d=|J|$ be the dimension of the cube $C(J,R)$. First we show that $\phi$ restricted to a cube $C(J,R)$ is a bijection. To this end, we need to introduce another concept and prove a lemma.

Let $i\in J$ and let $S_i=\{a_1\dts a_{2l}\}$. According to definition of the cube, there is an $h\leq l$ such that $S_i\setminus R=\{a_h,a_{h+l}\}$ and for every $X\in C(J,R)$, $|X\cap \{a_h,a_{h+l}\}|=1$. Hence there is a unique index $j$ such that
\bel{e-j}
\mbox{ either }\ \Pi_j(X)=\{a_h\dts a_{h+l-1}\},\
\mbox{ or }\ \Pi_j(X)=\{a_{h+1}\dts a_{h+l}\},
\ee
and this is an exclusive or. 
We will call such a $\Pi_j(X)$ the \emph{main projection in}~$S_i$.

\bll{l2.3}
For every $i\in J$, there is a unique $j$ such that for all $X\in C(J,R)$,  $\Pi_j(X)$ is the main projection in $S_i$.
\el
\bprf
First we observe that for every $i$ the cardinality of the set $X\cap(S_1\cup\dots\cup S_{i-1})$ does not depend on the choice of $X\in C(J,R)$. This is because, for every $i'\not\in J$, $X\cap S_{i'}=R\cap S_{i'}$ and for $i'\in J$, $|X\cap S_{i'}|=2l-1$. Hence the index $j$ of a projection that is a subset of $S_i$ does not depend on $X\cap(S_1\cup\dots\cup S_i)$. In particular, the indices $j$ of the projections $\Pi_j(X)$ that start with some $a_g$ with $g<h$ (we are using the notation above) only depend on $g$. This implies that the index of the main projection does not depend on which element of the cube we take.
\eprf

For a cube $C(J,R)$ and $i\in J$, we denote by $j(i)$ the index of the main projections in $S_i$.
Further, 
we denote the set of indices of the main projections by $B_{J,R}:=\{j(i)\ |\ i\in J\}$. In order to remember this concept, it is useful to note the duality: $J$ are indices of the sets of intervals of $S$ that are properly hit by some $X$ in the cube $C(J,R)$, while $B_{J,R}$ are indices of blocks of $X$ that properly hit some interval in $S$ (and these indices do not depend on a particular $X\in C(J,R)$).

\bll{l3.5}
Let $X,X'\in C(J,R)$, $X\neq X'$. Then $\phi(X)\neq\phi(X')$, in fact this still holds true if we restrict the vectors $\phi(X),\phi(X')$ to the coordinates in~$B_{J,R}$.
\el
\bprf
$X$ and $X'$ may be different only on some $S_i$ such that $i\in J$. For such an $i$, take $j=j(i)$. Then $\Pi_j(X)$ and $\Pi_j(X')$ are in a shift position. Since $\kappa$ is a coloring of the shift graph $Sh(N,l)$, we have $\kappa(X_j)\neq\kappa(X'_j)$, whence $\phi(X)\neq\phi(X')$.
\eprf

We will now estimate the number of different $\phi$ images of cubes.
Formally, we want to find an upper bound on the number of subsets of $[3]^p$ of the form
\bel{e-image}
\{(\kappa(X_1)\dts\kappa(X_p))\ |\ X\in C\}
\ee
for all cubes $C$ that can potentially appear in some ${S\choose k}$.
We will only need to bound the number of such sets for cubes of a sufficiently large dimension, but we would not save any essential factor if we restricted to such cubes.

Since each cube contains $2^d$ sets, where $d\leq k$ is the dimension of the cube with each element of the cube receiving one of the $3^p$ colors, there are, clearly, at most $3^{p2^d}$ possible images. However, we will need a better estimate for our purpose. While this trivial bound is doubly exponential in $k$, we will show a single exponential bound.

To count the number of images, we start by estimating the number of the sets $B_{J,R}=\{j_1\dts j_d\}\sub[p]$, which we can choose in ${p\choose d}$ ways. Next observe that for each $j_i$, $i=1,2\dts d$, there are precisely two sets $X_{j_i}$ when we range over all $X\in C(J,R)$. So we denote by $\{X_{j_i}^{(1)},X_{j_i}^{(2)}\}:=\{X_{j_i}\ |\ X\in C(J,R)\}$. Moreover, the $l$-sets $X_{j_i}^{(1)}$ and $X_{j_i}^{(2)}$ are in a shift position and, consequently, their $\kappa$ colors are different. This gives us
\bel{e(1)}
    {3\choose 2}^d=3^d
\ee
choices for the colors of $\{\kappa(X_{j_i})\ |\ i\in[d],k\in[2]\}$.

To address the case when $h\in[p]\setminus\{j_1\dts j_d\}$, let us first assume that $j_i<h<j_{i+1}$ for some $i=1,2\dts d-1$. We observe that once the values of the sets $X_{j_i}$ and $X_{j_{i+1}}$ are determined, so are also the values of sets
\bel{e-10}
X_{j_i+1},X_{j_i+2}\dts X_{j_{i+1}-1}.
\ee
The number of colorings of these sets is at most $3^{j_{i+1}-j_i}$ (in fact, it is smaller). Since $X_{j_i}\in\{X_{j_i}^{(1)},X_{j_i}^{(2)}\}$ and $X_{j_{i+1}}\in\{X_{j_{i+1}}^{(1)},X_{j_{i+1}}^{(2)}\}$  the pair $(X_{j_i},X_{j_{i+1}})$   attains four different values while $X$ is ranging over $X\in C(J,R)$. Thus we see only four different sequences~(\ref{e-10}) in the segment between $X_{j_i}$ and $X_{j_i+1}$, which have at most
\bel{e(2)}
(3^{j_{i+1}-j_i})^4
\ee
different colorings.

The cases of $h<j_1$ and $h>j_d$ are similar. Consequently we have at most
\bel{e(3)}
81^{p-d}
\ee
colorings of $(X_h;h\in[p]\setminus\{j_1\dts j_d\})$.

Combining (\ref{e(1)}),(\ref{e(2)}),\ref{e(3)} yields at most
\bel{e-bound}
\sum_{d=1}^n {p\choose d}\cdot 3^d\cdot 81^{p-d}<
\sum_{d=0}^p {p\choose d}\cdot 3^d\cdot 81^{p-d}
=(3+81)^p=84^p < 84^k,
\ee
different color patterns appearing on cubes.

One can, certainly, get a better bound (for instance, by using the fact that all consecutive projections $X_j,X_{j+1}$ are in a shift position, hence their colors must be different), but we will satisfy ourselves with this bound, because the better one would have very little influence on the bound that we will get.


\subsection{The rest of the proof of Theorem \ref{t-x1}}

Since the image $\phi(C)$ of a cube $C$ of dimension $d$ has $2^d$ elements, Chernoff's bound gives us that the probability that a random $\psi:[2^d]\to[2]$ has relative discrepancy $>\delta$ on $C$ is bounded by $<2\exp(-2\delta^22^d)$.
We are only concerned with cubes of dimensions $\geq n/2\e$.
Setting  $d=n/2e$,
we get the following bound 
\[
<2\exp(-2\delta^22^{(n/2\e)}) = 2\exp(-2\delta^22^{(k/2\e(2l-1))})<\exp(-2\delta^22^{k/12 l}).
\]
By the union bound, if
\bel{e-union}
84^k\cdot \exp(-2\delta^22^{k/12 l})=\exp(-2\delta^22^{k/12 l}+k\ln 84)<1,
\ee
then there exists a mapping $\psi$ that has relative discrepancy $<\delta$ on all images of cubes of dimensions $\geq n/2\e$.

Since $\ln 84<5$, in order to prove~(\ref{e-union}), it will suffice to show that
\[
5k-2\delta^22^{k/12 l}<0.
\]
Setting $\delta=2^{-k/60 l}$ this translates to
\[
5k<2^{\frac k{12 l}-\frac k{30 l}}=2^{\frac k{20 l}},
\]
or
\bel{e-l-bound}
l<\frac 1{20}\frac k{\log_2 5k},
\ee
which holds true in view of our assumption in Theorem~\ref{t-x1}.
For $k$ sufficiently large, we can replace the right hand side with $\frac 1{21}\frac k{\log k}$. 
So we proved:
\bfa
If the shift graph $Sh(N,l)$ can be colored by 3 colors and  $l<\frac 1{21}k/\log_2 k$, where $k$ and $l$ are sufficiently large, then there exists a coloring $\psi$ with discrepancy $\leq 2^{-k/60 l}$ on every cube of dimension at least~$n/2\e$.
\efa
By Lemma~\ref{l-shift}, a 3-coloring of $Sh(N,l)$  is guaranteed if $N<4\tw_{l-4}(2)$. Also Corollary~\ref{c3.1} yields that all but a fraction of
\[
\exp(-n/30)=\exp(-\frac k{60l-30})<\exp(-\frac k{60l})
\]
$k$-subsets of $S$ are covered by cubes of dimensions $\geq n/2\e$. Thus, by Fact~\ref{fact}.3, the relative discrepancy of $\gamma=\psi\circ\phi$ on ${S\choose k}$ is at most
\[
\exp(-\frac k{60l})+\delta=\exp(-\frac k{60l})+2^{-k/60l}<2^{-k/61l},
\]
where the inequality holds if $k/l$ is sufficiently large, but due to~(\ref{e-l-bound}), $k/l\mapsto \infty$ when $k\mapsto\infty$. Hence it suffices to assume that $k$ is sufficiently large. 
\qed


\subsection{Proof of Theorem \ref{t-1} from Theorem \ref{t-x1}}

We will need the following fact.
\bfa\label{f-1}
Let $A,B$ be disjoint sets with $|A|=k+n$, $|B|=2n$, where $n\leq k/4$. Then  the cardinality of the set
\[
{\cal X}:=\{X\in{A\cup B\choose k}\ |\ |B\cap X|<n\}
\]
is less than
\[
\alpha^n{k+3n\choose k},
\]
where $\alpha\leq\exp(-\frac 2{147})$.
\efa
\bprf
The bound follows from the estimate on the bound on the hypergeometric distribution of Theorem~\ref{t-hyper}. The population size is $L=k+3n$, the number of success states is $2n$ and the number of draws is~$k$. Plugging the numbers in (\ref{e-hyper}), with $p=\frac{2n}{k+3n}$ and
\[
\tau=p-\frac nk=\frac{n(k-3n)}{k(3n+k)},
\]
we obtain
\[
|{\cal X}|\leq \exp(-2\tau^2\frac{k^2}{L-k})=
\exp\left(-\frac 23 n(\frac{3n-k}{3n+k})^2\right).
\]
Assuming $n\leq k/4$, this can be upper-bounded by $\exp(-\frac 2{147}n)$.
\eprf

Recalling that unlike in Theorem~\ref{t-x1}, Theorem~\ref{t-1} has no
divisibility assumption. So our next goal is to remove this
assumption. For a given $n$ satisfying the assumptions of Theorem~\ref{t-1}, let $k'$ be the largest odd multiple of $n$ less than or equal to $k$. Since $k\geq 11n$ by an assumption of Theorem~\ref{t-1} and $k-k'<2n$, we conclude that $k'=(2l-1)n\geq 9n$ and hence also $l\geq 5$. We intend to apply Theorem~\ref{t-x1} with $k'$. Note that for a sufficiently large constant $c$, the condition $c\log k\leq n$ of Theorem~\ref{t-1} implies the condition $k'/l\geq 165\ln 2k'$ of Theorem~\ref{t-x1}. Hence, if $N\leq\tw_{l-4}(2)$, then there exists a mapping $\gamma':{[N]\choose k'}\to\{-1,+1\}$ such that for every set $S\sub[N]$ of size $k'(1+\frac 1{2l-1})$, we have  $disc(\gamma',{S\choose k'})\leq 2^{-k'/61l}$. We define $\gamma:{[N]\choose k'}\to\{-1,+1\}$ by
\[
\gamma(x_1,x_2\dts x_k)=\gamma'(x_1,x_2\dts x_{k'})
\]
for every $k$-tuple $1\leq x_1<x_2<\dots<x_{k'}<\dots<x_k\leq N$. We will show that for every $M\sub [N]$, $|M|=m=k+3n$, $disc(\gamma,{M\choose k})\leq 2^{-\eps n}$ for $\eps=1/213$. (For Theorem~\ref{t-1}, we need $|M|=k+n$ and also $\lfloor k/n\rfloor$ instead of~$l$; we will address this issue later.)

Let $M=\{y_1<y_2<\dots<y_m\}$, $A:=\{y_1<y_2<\dots<y_{k+n}\}$ and  $B=\{y_{k+n+1}<y_{k+n+2}<\dots<y_m\}$ where $y_{k+n}<y_{k+n+1}$. Let
\[
{\cal X}_{bad}:=\{X\in {M\choose k}\ |\ |B\cap X|<n\}.
\]
By Fact \ref{f-1}, $|{\cal X}_{bad}|\leq \alpha^n{m\choose k}$ with 
$\alpha=\exp(-\frac 2{147})$.
For $X\in{M\choose k}\setminus {\cal X}_{bad}$, let $Q_X$ be the set of its last $k-k'$ elements and observe that for each such $X$, $Q_X\sub B$. Next we divide all $X\not\in{\cal X}_{bad}$ into groups sharing the last $k-k'$ elements. For each $Q\in{B\choose k-k'}$, let
\begin{align*}
{\cal X}_Q&:=\{X\in{M\choose k} |\ Q_X=Q\},\\
S_Q&:=\{y\in M\ |\ y<\min Q\}.
\end{align*}
Since $X\not\in{\cal X}_{bad}$, we have $|S_{Q_X}|> k+n$.
Therefore, for every such $Q$,
\[
disc(\gamma,{\cal X}_Q)=disc(\gamma',{S_Q\choose k'})\leq 2^{-k'/61l}.
\]
Hence
\[
disc(\gamma,{M\choose k})\leq  2^{-k'/61l} +\frac{|{\cal X}_{bad}|}{{m\choose k}}
\leq  2^{-k'/61l} +\exp(-\frac{k'}{147l})<2^{-k'/213l}.
\]

Summarizing, for $5\leq l\leq \frac 1{165}\frac k{\log 2k}$, $k\geq \max\{k_0,11n\}$ and $N\leq \tw_{l-4}(2)$, we just proved that there exists a 2-coloring $\gamma$ with $disc(\gamma,{M\choose k})\leq 2^{-k/213l}$ for every $M\sub[N]$, $|M|=k+3n$. Since
\(
\frac kl\geq\frac k{2l-1}\geq \frac{k'}{2l-1}=n,
\)
we can bound the discrepancy by $2^{-\frac n{213}}$.

Now we set $n':=3n$ and recalling that $k'=(2l-1)n$ we observe that
\[
\frac k{n'}=\frac k{3n}\leq \frac{k'+2n}{3n}\leq\frac{k'+n}{2n}-4=l-4,
\]
which implies that mapping $\gamma$ with the property stated above validates the statement of Theorem~\ref{t-1}.




\section{An explicit construction}\label{s-explicit}

\bt\label{t-2}
For every odd $c\geq 3$, there exist $l_0$ 
such that for every $k,N$, $l\geq l_0$ and $0<\delta<1$ such that $k$ is divisible by $l$, there exists $\eps>0$ such that the following holds true. If
$$
k\geq l^{7+3\delta}\quad\mbox{ and }\quad N\leq 4\tw_{l-4}(2),
$$
then one can explicitly define a coloring $\rho$ of $k$-element subsets of $[N]$ by $c$ colors such that  for every subset $S\sub [N]$ of size $|S|\geq k(1+1/l)$, the relative discrepancy of $\rho$ restricted to ${S\choose k}$ is
\[
disc(\rho,\mbox{${S\choose k}$})\leq\exp(-\eps l^{\delta}).
\]
\et
The construction is explicit in the following sense. The coloring is constructed from two mappings. The first one is a vector of colors in a $3$-coloring of $Sh(N,l)$. The $3$-coloring $\chi$ of $Sh(N,l)$ in Lemma~\ref{l-shift} is constructed recursively starting with the trivial coloring of  $Sh(N,1)$ by $N$ colors and using a simple construction that produces a coloring of $Sh(N,h+1)$ from a coloring of  $Sh(N,h)$, see Section~\ref{s-shift}. While we do not attempt to give a formal definition, we believe  it is reasonable to call such a construction explicit.
If we agree that the 3-coloring is explicit, then also the first mapping is explicit.
The second mapping is clearly explicit because it is just addition modulo~$c$ of the elements of the vector produced by the first coloring.

If we disregard the constructivity, then Theorem~\ref{t-2} is weaker than Theorem~\ref{t-1}, because $l$ has to be much smaller than $k$, consequently $N$ also has to be smaller. However it should be noted that we have not striven for the best bound. We prefer to give a simpler construction rather than the best possible, one reason being that we do not see how to get the same bound as in Theorem~\ref{t-1} anyway.

To prove Theorem~\ref{t-2}, we will need the following lemma.

\bll{l-parity}
Let $0<p<1$, $q=1-p$, let $l,h$ be integers $0\leq h<l$, and let $\bar x=x_1\dts x_n$ be independent Bernoulli variables with $\prob[x_i=1]=p$ and $\prob[x_i=0]=q$. Then
\[
\left|\prob[\ \sum_i x_i\equiv h\mod l\ ] - \frac 1l\right|\ \leq\ (1-2pq(1-\cos(2\pi/l)))^{n/2}.
\]
\el
\bprf
Let $\beta:=2\pi/l$ and let $\alpha:=\cos\beta+\iu\sin\beta$, one of the two primitive $l$-th roots of unity closest to~1. Then
\begin{align*}
\sum_{j=0}^{l-1}\alpha^{(i-h)j}&=0\mbox{ if }i\not\equiv h \mod l\\
                           & =l\mbox{ otherwise. }
\end{align*}
Thus
\begin{align*}
\prob[\ \sum_i x_i\equiv h\mod l\ ]
&=\sum_{i=0}^n\binom{n}i p^iq^{n-i}\frac 1l\sum_{j=0}^{l-1}\alpha^{(i-h)j}\\
&=\frac 1l\sum_{j=0}^{l-1}\alpha^{-hj}\sum_{i=0}^n\binom{n}i \alpha^{ij} p^iq^{n-i}\\
&=\frac 1l\sum_{j=0}^{l-1}\alpha^{-hj}(\alpha^jp+q)^n\\
&=\frac 1l(1+\sum_{j=1}^{l-1}\alpha^{-hj}(\alpha^jp+q)^n)
\end{align*}
The largest terms in absolute values in the sum are for $j=1$ and $j=l-1$ and their absolute values are $|\alpha p+q|^n$. Thus it suffices to compute $|\alpha p +q|$. 
\begin{align*}
|\alpha p+q|& =|q+p\cos\beta+\iu p\sin\beta|\\
  &=\sqrt{(q+p \cos\beta)^2+(p \sin\beta)^2}\\
&=\sqrt{q^2+p^2+2pq\cos\beta}\\
&=\sqrt{1-2pq(1-\cos\beta)}
\end{align*}
\eprf

\paragraph{Set-up of parameters.} Let $k$ and $l$ be given. For the rest of the proof of Theorem~\ref{t-2} we set 
\[
n:=k/l,\quad m:=n(l+1).
\]
Given an $S\sub[N]$, $|S|=m$, we divide $S$ into $n$ segments $S_i$ of length $l+1$. In this part we will use a different definition of properly hitting.
\bdf
For $X\sub S$, we say that $X$ \emph{properly $h$-hits} $S_i$ if
\ben
\item (parity condition)  $|X\cap (S_1\cup\dots\cup S_{i-1})|=h\mod l$, and
\item (hitting condition) $X\cap S_i$ is either $S_i\setminus\{\max S_i\}$, or $S_i\setminus\{\min S_i\}$.
\een
We say that $X$ \emph{properly hits} $S_i$ if $X$ {properly $0$-hits} $S_i$.
\edf
Sometimes, as in the next lemma, it may be easier to work with 0/1-strings, so 
if $\bar x\in\{0,1\}^m$, we say that $x$ \emph{$h$-properly hits} $S_i$ if the set of coordinates on which $\bar x$ is equal to~1  $h$-properly hits $S_i$.


The reason for introducing the concept of hitting with the parameter~$h$ is that we will need to get a bound for Azuma'a inequality. To this end we will, for every $1\leq t\leq n$, divide the sequence of intervals $S_i$ into two parts
\[
S_1\dts S_{t}\mbox{ and }S_{t+1}\dts S_n
\]
and estimate the expected number of properly hit intervals in the second part assuming the $X\cap(S_1\cup\dots\cup S_{t})$ is given. Clearly, for $i> t$, $S_i$ is properly hit iff $S_i$ is $h$-properly hit when we consider only the part $S_{t+1}\dts S_n$ and $h=-|X\cap(S_1\cup\dots\cup S_{t})|\mod l$.

So we need to estimate the expected number of $h$-properly hit intervals in $S_{t+1}\dts S_n$. In the following lemma, in order to simplify notation, we will not index intervals $S_i$ from $t+1$ to $n$, but from $1$ to~$n'$, where $n'=n-t$.


\bll{l-4.2}
For every $\delta>0$, there exists $l_0$ such that the following holds true for all $l\geq l_0$, $n'$ and $0\leq h<l$. Let $m'=n'(l+1)$ and 
let  $\bar x=x_1\dts x_{m'}$ be a string of independent variables with $p:=\prob[x_i=1]=l/(l+1)$ and $q:=\prob[x_i=0]=1/(l+1)$, and let $Z_{n',h}$ be the random variable counting the number of segments $S_i$ properly $h$-hit by $\bar x$. Then
\[
(2{\e}^{-1}-o(1))\left(\frac{n'}{l^2}-l^{\delta}\right)\leq
\mathbb{E}(Z_{n',h})
\leq (2{\e}^{-1}+o(1))\left(\frac{n'}{l^2}-l^{\delta}+l^{1+\delta}\right),
\]
where $o(1)$ only depends on $l$.%
\el
\bprf
The condition that $\bar x$ $h$-properly hits $S_i$ is a conjunction of two independent events: $A_i$, the event that $\bar x$ satisfies the parity condition for $S_i$, and $B_i$, the event that $\bar x$ satisfies the hitting condition for $S_i$. 

To compute $Pr[A_i]$ we will bound 
$(1-\cos\frac{2\pi}l)\leq\frac{2\pi^2}{3l^2}$ (assuming $l$ is sufficiently large) and apply Lemma~\ref{l-parity} to infer that

\begin{align*}
\left|\prob(A_i)-\frac 1l\right|\leq & \left(1-\frac{2l}{(l+1)^2}\cdot\frac{2\pi^2}{3l^2}\right)^{i(l+1)/2}\\
<&\exp\left(-\frac{2\pi^2i}{3(l+1)l}\right)\\
<&\exp\left(-\frac{\pi^2i}{2l^2}\right).
\end{align*}
Consequently, for $i>l^{2+\delta}$, we have
\bel{e-star}
\left|\prob(A_i)-\frac 1l\right|<\exp(-2\pi^2l^\delta)=o(1).
\ee
On the other hand
\[
\prob(B_i)=\frac 2{l+1}\left(\frac l{l+1}\right)^l=(1+o(1))\frac{2}{(l+1)\e}> \frac{2}{(l+1)\e}.
\]
In order to bound $\mathbb{E}(Z_{n',h})$, we first observe that
\[
\mathbb{E}(Z_{n',h})=\sum_{i=1}^{n'}\prob(A_i\wedge B_i)=\sum_{i=1}^{n'}\prob(A_i)\prob(B_i).
\]
To prove the lower bound, we ignore intervals $S_i$ for $1\leq i< l^{2+\delta}$. Beyond this segment, we use~(\ref{e-star}) to infer that
\begin{align*}
\mathbb{E}(Z_{n',h})&\geq \sum_{i> l^{2+\delta}}^{n'}\prob(A_i)\prob(B_i)\\
& \geq\sum_{i> l^{2+\delta}}^{n'} \left(\frac 1l-o(1)\right)\frac{2}{\e (l+1)}\\
& = (1-o(1))\cdot (n'-l^{2+\delta})\frac{2}{\e l^2}\\
& = (1-o(1))\cdot 2\e^{-1}\left(\frac {n'}{l^2}-l^\delta\right).
\end{align*}
For the upper bound, we ignore the parity condition up to $i= \lfloor l^{2+\delta}\rfloor$. Thus we get
\begin{align*}
\mathbb{E}(Z_{n',h})&\leq \sum_{i\leq l^{2+\delta}}\prob(B_i)   +\sum_{i> l^{2+\delta}}^{n'}\prob(A_i)\prob(B_i)\\
&\leq (1+o(1))\cdot\frac{2}{\e (l+1)}l^{2+\delta}+(1+o(1))\cdot 2\e^{-1}(\frac{n'}{l^2}-l^\delta)\\
& \leq (1+o(1))\cdot 2\e^{-1}\left(l^{1+\delta}+\frac{n'}{l^2}-l^{\delta}\right).
\end{align*}
\eprf

Observe that the bounds do not depend on~$h$ and the difference between the upper and lower bounds is 
\bel{e-regard}
(1+o(1))2\e^{-1}l^{1+\delta}
\ee
and it does not depend on $n'$. We will use it for Azuma's inequality.

Write $(x_1,x_2\dts x_m)$ as $(y_1,y_2\dts y_n)$, where $y_i$ are segments of $l$ consecutive elements of $(x_1,x_2\dts x_m)$. 
Denote by $f(y_1\dts y_n)$ the number of properly hit segments $S_j$. We need to find an upper bound on
\bel{e-diff}
| \mathbb{E}[f(y_1\dts y_{t-1},y_t,Y_{t+1}\dts Y_n)]-  \mathbb{E}[f(y_1\dts y_{t-1},y'_t,Y_{t+1}\dts Y_n)]|
\ee
for every $t$, $1\leq t\leq n$. Let $t$ be given. Then
\[
\mathbb{E}[f(y_1\dts y_{t-1},y_t,Y_{t+1}\dts Y_n)]=A+B,
\]
where $A$ is the number of properly hit intervals $S_j$ for $j=1\dts t$ and $B$ is $\mathbb{E}(Z_{n-t,h})$, where $h=-\sum_{t=1}^{t(l+1)}x_t \mod l$. Similarly decompose 
\[
\mathbb{E}[f(y_1\dts y_{t-1},y'_t,Y_{t+1}\dts Y_n)]=A'+B'.
\]
Clearly $|A-A'|\leq 1$, because $S_t$ is the only interval on which one of the strings may properly hit it and the other not. $B'$ may be different from $B$ if $h'=-(\sum_{t=1}^{t(l+1)-1}x_t)-x'_t\not\equiv h\mod l$, but the difference is at most the difference between the lower bound and upper bounds on~$\mathbb{E}(Z_{n-t,h})$, which is~(\ref{e-regard}). Thus (\ref{e-diff}) is bounded by
\bel{e-K}
1+(1+o(1))2\e^{-1}l^{1+\delta}=(1+o(1))2\e^{-1}l^{1+\delta}.
\ee
Hence, in the following lemma, we can apply Azuma's inequality~(\ref{e-azuma}) with $K=(1+o(1))2\e^{-1}l^{1+\delta}$.

\bll{l-4}
Let $n\geq l^{6+3\delta}$, with $l$ sufficiently large. Then with probability
\[>1- 2\exp(-(1-o(1))l^\delta/8)\]
the number of segments properly hit is at least $n/\e l^2$.
\el
\bprf
For $h=0$ and $n'=n$, Lemma~\ref{l-4.2} gives us an estimate on the expected value of properly hit intervals. With our assumption $n\geq l^{6+3\delta}$, the terms $l^\delta$ and $l^{1+\delta}$ become negligible, so $\lambda=n/\e l^2$ is approximately $1/2$ of the expected value. When we apply Azuma's inequality (Theorem~\ref{l-azuma}) with this $\lambda$ and $K$ estimated above~(\ref{e-K}), the exponent in the inequality becomes (omitting factors $(1-o(1))$)
\[
\frac{\lambda^2}{2K^2n}=\frac{n^2}{2\cdot\e^2l^4\cdot 4\e^{-2}l^{2+2\delta}n}=\frac n{8l^{6+2\delta}}\geq \frac{l^\delta}{8}.
\]
\eprf

We now consider random subsets of $S$ of size~$k$. 
In the same way as in Lemma~\ref{l-number-of-hit}, we get from Lemma~\ref{l-4}:

\bll{l-random-subsets}
Let $n\geq l^{6+3\delta}$. For $X$, a random subset of $S$ size $k$, $X$ properly hits at least $n/\e l^2$ segments $S_i$ with probability $1-2\exp(-(1-o(1)) l^\delta/8)$.
\el

We now define cubes similarly as in Definition~\ref{d-cube}.
\bdf\label{d-cube2}
Let $J\sub[n]$ be nonempty and $R\sub S$ such that
\bi
\item for every $j\in J$, $R\cap S_j=\emptyset$ and  $|R\cap (S_1\cup\dots\cup S_{j-1})|=0\mod l$,
\item $|R|+l\cdot|J|=k$.
\ei
Then \emph{the cube $C'(J,R)$} is the set of $k$-element subsets $X\sub S$ such that
\bi
\item for every $j\in J$, $X$ properly hits $S_j$,
\item for all $j\in[n]\setminus J$, $X\cap S_j=R\cap S_j$.
\ei
The {dimension of the cube} is $|J|$.
\edf
If $X$ properly hits $d$ intervals $S_i$, then, clearly, there is a maximal cube $C'(J,R)$ of dimension~$d$ that contains~$X$. 
Thus Lemma~\ref{l-random-subsets} guarantees that ${S\choose k}$ is covered by $k$-sets of disjoint cubes of dimensions $\geq n/\e l^2$ except for an exponentially small fraction of the elements of ${S\choose k}$, as stated in the following corollary.

\begin{corollary}\label{cor4.1}
All $k$-sets in ${S\choose k}$, except for a fraction of
$$
2\exp(-(1-o(1)) l^\delta/8),
$$
are contained in disjoint cubes of dimensions  $\geq n/\e l^2$.
\end{corollary}

Let $\kappa$ be a proper coloring of the shift graph $Sh(N,l)$ by $3$~colors; we now prefer to use colors $0,1,2$. The coloring $\phi$ of ${N\choose k}$ by $3^n$ colors is defined as follows:
\[
\phi(X)=(\kappa(X_1)\dts \kappa(X_n)),
\]
where $X_1\dts X_n$ is a partition of $X$ into $n$ consecutive segments of length~$l$. (Recall that $k=nl$.)

The key difference between the non-explicit and explicit constructions is that we define the second coloring in a different way---we simply add the colors in $\phi(X)$ modulo~$c$. Thus the coloring $\rho$ of $[N]\choose k$ is defined by
\[
\rho(X):= \sum_{i=1}^n\kappa(X_i) \mod c.
\]
Since Lemma \ref{l-random-subsets} guarantees that all $k$-element sets, except for an exponentially small fraction of them, are covered by disjoint cubes of dimensions $\geq n/\e l^2$, it suffices to show that $\rho$ has exponentially low discrepancy on every cube of dimension~$\geq n/\e l^2$.

Let $C$ be a cube of dimension $d\geq n/\e l^2$ and let $S_i$ be a segment properly hit by the sets of $C$. For a set $X\in C$, one of the segments $X_1\dts X_n$ properly hits $S_i$; let this be $X_j$. As in the previous construction, $j$ is uniquely determined by $i$. As before, let $B_C$ denote the set of such indices~$j$, i.e.,
\[
B_C:=\{j\in[n]\ |\ \exists i\in[n]\ \forall X\in C,\ X_i\mbox{ properly hits } S_j\}. 
\]
For each $j\in B_C$, there are two possibilities for $X_j$. The two possible versions are in a shift position, so they have two different colors $a_j,b_j\in\{0,1,2\}$. We consider two possibilities:
\ben
\item $|a_j-b_j|=1$,
\item $|a_j-b_j|=2$.
\een
According to these possibilities, we split $B_C$ into two parts $B^1_C$ and $B^2_C$. Let $d_1:=|B^1_C|$ and $d_2:=|B^2_C|$, so $d_1+d_2=d$. Consider two cases.

Case 1, $d_1\geq d/2$. Then we split cube $C$ into disjoint cubes of dimension $d_1$ by fixing $X_j$s on every $j\in B^2_C$. Let $C'$ be one of these $2^{d_2}$ cubes. Observe that the elements of $C'$ are in 1--1 correspondence with sequences $(\xi_1,\xi_2\dts \xi_{d_1})\in\{0,1\}^{d_1}$ with each $X\in C'$ being assigned $\xi_j\in\{0,1\}$ depending on the color of the $l$-set~$X_j$. Consequently, the colors $\rho(X)$ for $X\in C'$ can be represented by
\[
A+\sum_{j=1}^{d_1}\xi_j \mod c,
\]
where $A$ is constant and $\xi_j\in\{0,1\}$. According to Lemma~\ref{l-parity}, there exists a constant $\alpha_c:=1-\frac 12(1-\cos(2\pi/c))<1$ such that the frequency of every modulus is $1/c$ up to an error $\alpha_c^{d_1}$. Since  $d_1\geq n/2\e l^2> l$, this error is negligible when compared with our target error $\exp(-\eps l^{\delta})$. So, like in the first construction, the main error comes from the fraction of $k$-sets that are not covered by large cubes, which is estimated in Corollary~\ref{cor4.1}.

Case 2, $d_2\geq d/2$. The same argument shows that the colors $\rho(X)$ in a cube of dimension $d_2$ can be represented by
\[
A+\sum_{j=1}^{d_2}2\xi_j \mod c,
\]
but because in $\Z_c$ we can divide by 2, we get the same conclusion.

This concludes the proof of Theorem~\ref{t-2}.
\qed


\section{Proof of Lemma~\ref{l-shift}.}\label{s-shift}

Recall the lemma:
\bl[\ref{l-shift} restated]\label{l5.1}
Let $l\leq N\leq 4\tw_{l-4}(2)$. Then the chromatic number of $Sh(N,l)$ is~$\leq 3$ and, for $N\geq 2l+1$, we get equality. Moreover, such a $3$-coloring can be explicitly defined.
\el

We start by introducing some notation and making preliminary considerations.
A proper coloring of a directed graph $G$ is a proper coloring of $G$ with the orientation of the edges ignored. We will denote by $\chi(G)$ the chromatic number of $G$, i.e., the least number of colors needed to color the graph. There is a natural orientation of the edges of a shift graph which makes it a directed acyclic graph. We will use the same notation for the directed shift graphs. 
Although directed shift graphs are acyclic, we will also need graphs with cycles. 


For a directed graph $G$, we denote by $\partial G$ the graph whose vertices are edges of $G$ and the edges of $\partial G$ are oriented paths of length~2.

We note that  for $l\geq 1$, $Sh(N,l)\cong \partial^{l-1}Sh(N,1)$. For example, the vertices of $\partial^2Sh(N,1)$ are $((a,b),(b,c)$), where $1\leq a<b<c\leq N$, and they correspond to vertices $\{a,b,c\}$ in $Sh(N,3)$.


Furthermore, we will need the tower-of-square-roots-of-two function, which is defined by
\[
\overline{\tw}_1(x):=2x\quad\mbox{ and }\quad\overline{\tw}_{i+1}:= 
2^{\frac 12\overline{\tw}_i(x)}=\sqrt 2^{\overline{\tw}_i(x)}.
\]
We observe that $\overline{\tw}_i(x)$ is an even number for all integers $i\geq 0$ and $x\geq1$.

Below we summarize some well-known results and, for completeness, we also  sketch their proofs. Clearly $\chi(Sh(N,1))=N$.
The chromatic number of shift graphs for larger~$l$ can be estimated using the following well-known recurrence relations~\cite{erdos-hajnal-A,harner-entringer,poljak-rodl}.

\bll{l5.2}
Let $s\geq 1$ be an integer. Then 
\ben
\item $\chi(G)\leq 2^s$ implies $\chi(\partial G)\leq 2s$;
\item $\chi(Sh(N,h))\leq {2s\choose s}$ implies $\chi(Sh(N,h+1))\leq 2s$.
\een
\el
\bprf
1. If $\chi(G)\leq 2^s$, then the edge set can be covered by $s$ bipartite graphs. Each of these bipartite graphs can be further decomposed into two bipartite graphs in which each edge has the tail in one set of the bipartition and the head in the complement. Such sets are independent in $\partial G$. Hence $\chi(\partial G)\leq 2s$.

\medskip
2. Let $\alpha$ be a coloring of  $Sh(N,h)$ by ${2s\choose s}$ colors. Suppose that the colors are subsets of size $s$ of~$[2s]$. We define a coloring $\beta$ of $Sh(N,h+1)$ by
\[
\beta(\{x_1,x_2\dts x_h,x_{h+1}\})=c,
\]
where $c\in[2s]$ is an arbitrary (say, the first) element such that
\[
c\in \alpha(\{x_1,x_2\dts x_h\})\setminus\alpha(\{x_2\dts x_h,x_{h+1}\}).
\]
We leave to the reader to verify that this defines a proper coloring of ${Sh(N,h+1)}$.
\eprf

Since $\partial(Sh(N,h))\cong Sh(N,h+1)$, Lemma~\ref{l5.2}.1 implies
\bel{e-implication}
\mbox{If }\chi(Sh(N,h))\leq\overline{\tw}_{i+1}(s)\mbox{ then }
\chi(Sh(N,h+1))\leq\overline{\tw}_{i}(s).
\ee
(Recall that $\overline{\tw}_{i+1}(s)=2^{\frac 12\overline{\tw}_{i}(s)}$ and $\overline{\tw}_{i}(s)$ is an even number.)
Thus $(l-2)$-times repeated application of the above implication (\ref{e-implication}) yields that
\[
N=\chi(Sh(N,1))\leq\overline{\tw}_{l-3}(s)
\mbox{ implies }
\chi(Sh(N,l-3))\leq\overline{\tw}_{1}(s)=2s.
\]
For $s=4$, this gives us:
\bfa
If $N\leq \overline{\tw}_{l-4}(4)$, then $\chi(Sh(N,l-4))\leq 8$.
\efa
Since $\chi(Sh(N,l-3))\leq 2^3$,
applying Lemma~\ref{l5.2}.1 again yields $\chi(Sh(N,l-2))\leq 6$. Now we apply Lemma~\ref{l5.2}.2 and get $\chi(Sh(N,l-1))\leq 4$. Given a $4$-coloring of $Sh(N,l-1)$, we get a homomorphism $\phi$ of $Sh(N,l-1)$ into $\vec K_4$, which is the directed graph on four vertices with all twelve ordered pairs as edges. As is well-known (see~\cite{poljak,schmerl,duffus-lefmann-rodl}), $\partial^2 \vec K_4$ can be colored by $3$ colors. Such a coloring $\psi$ is defined by
\begin{align*}
\psi((a,b),(b,c)):&= b\mbox{ if }b\neq 4\\
&= d\in \{1,2,3\}\setminus\{a,c\}\mbox{ if }b=4.
\end{align*}
Thus composing $\phi$ with $\psi$ we get a $3$-coloring of $Sh(N,l)$, which proves

\bfa\label{f5.2}
If $N\leq \overline{\tw}_{l-4}(4)$, then $\chi(Sh(N,l))\leq 3$.
\efa
If $N\geq 2l+1$, then $Sh(N,l)>2$ because $Sh(2l+1,l)\sub Sh(N,l)$ contains an undirected cycle of length $2l+1$.

To finish the proof of Lemma~\ref{l5.1}, it remains to estimate from below the tower function $\overline{\tw}$ using the standard tower function~$\tw$. Toward that end we observe that for $x\geq 2$,
\bel{e-xy}
2^{2x}=\sqrt 2^{4x}\geq 4\cdot 2^x=2^{x+2}.
\ee
This is equivalent to $\overline{\tw}_2(2x)\geq 4\tw_2(x)$. From the above inequality, we get by induction
\[
\overline{\tw}_{i+1}(2x)=\sqrt{2}^{\overline{\tw}_i(2x)}
\geq\sqrt{2}^{4\cdot{\tw}_i(x)}\geq 4\cdot 2^{\tw_i(x)}=4\cdot\tw_{i+1}(x),
\]
where the first inequality is the induction assumption and the second one follows from~(\ref{e-xy}). 
Hence we have, for all $i\geq 1$,
\[
\overline{\tw}_i(2x)\geq 4\tw_i(x).
\]
In particular, for $x=2$ and $i=l-4$, we get $\overline{\tw}_{l-4}(4)\geq 4\tw_{l-4}(2)$. This inequality and Fact~\ref{f5.2} immediately yield the bound stated in Lemma~\ref{l5.1}.

Since the colorings in Lemma~\ref{l5.2} can be explicitly defined, the 3-coloring of $Sh(N,l)$ can also be explicitly defined.
\qed









\section{Conclusions and open problems}

The natural question is how much the bound on $N$ can be relaxed. This would be interesting even if the bound on the subsets $S$ is also relaxed. 
In particular we would like to know the answer to the following problem. (We believe that the answer is positive.)

\begin{problem}\label{p-1}
Does there exist $\eps>0$ such that for every $k$ sufficiently large and $ N\leq \tw_{\lfloor \eps k\rfloor}(2),$ there exist a $2$-coloring $\gamma$ of ${[N]\choose k}$  such that  for every subset $S\sub [N]$ of size $|S|=2k$, $disc(\gamma,{S\choose k})=o(1)$?
\end{problem}

\end{document}